\documentclass[11pt,preprint]{imsart}
\RequirePackage[OT1]{fontenc}
\RequirePackage{amsthm,amsmath}

\usepackage{graphicx,amsmath,amsthm,amsfonts,amssymb}
\usepackage{booktabs, rotating, float}
\usepackage[hmargin=3cm,vmargin=3cm]{geometry}
\usepackage{mathtools}
\usepackage{bbm}
\mathtoolsset{showonlyrefs}
\usepackage{mathrsfs}
\usepackage[authoryear]{natbib}
\usepackage{enumerate}
\usepackage{color}
\usepackage{listings}
\usepackage{epstopdf}
\usepackage{makecell}
\usepackage{tikz}
\usetikzlibrary{decorations.pathreplacing}
\usepackage{array}
\usepackage{caption}
\usepackage{subcaption}
\usepackage{soul}
\usepackage{hyperref}
\usepackage{verbatim}
\usepackage{siunitx}
\usepackage{graphicx}
\usepackage[font=scriptsize]{caption}
\usetikzlibrary{calc} \tikzset{>=latex}
\usetikzlibrary{calc,decorations.markings}
\usetikzlibrary{positioning,shapes,decorations.text,patterns}
\theoremstyle{plain}
\newtheorem{thm}{Theorem}[section]
\newtheorem{cor}[thm]{Corollary}

\newtheorem{prop}[thm]{Proposition}
\newtheorem{lemma}[thm]{Lemma}
\newtheorem{assumption}[thm]{Assumption}
\theoremstyle{definition}
\newtheorem{defn}[thm]{Definition}
\newtheorem{rem}[thm]{Remark}
\newtheorem{eg}[thm]{Example}
\newtheorem{fact}[thm]{Fact}
\newtheorem{observe}[thm]{Observation}
\newtheorem{conj}[thm]{Conjecture}

\numberwithin{equation}{section}

\newcommand{\rpm}{\sbox0{$1$}\sbox2{$\scriptstyle\pm$}
  \raise\dimexpr(\ht0-\ht2)/2\relax\box2 }

\tikzstyle{nd} = [anchor=base, inner sep=0pt]
\tikzstyle{ndpic} = [remember picture, baseline, every node/.style={nd}]

\def\beq{\begin{equation}}
\def\eeq{\end{equation}}
\def\ba{\begin{enumerate}[(a)]}
\def\bei{\begin{enumerate}[(i)]}
\def\be{\begin{enumerate}[(1)]}
\def\ee{\end{enumerate}}
\def\bi{\begin{itemize}}
\def\ei{\end{itemize}}
\def\beg{\begin{eg}}
\def\eeg{\end{eg}}
\def\bd{\begin{defn}}
\def\ed{\end{defn}}
\def\bt{\begin{thm}}
\def\et{\end{thm}}
\def\bl{\begin{lemma}}
\def\el{\end{lemma}}
\def\bfac{\begin{fact}}
\def\efac{\end{fact}}

\def\bc{\begin{cor}}
\def\ec{\end{cor}}
\def\bp{\begin{prop}}
\def\ep{\end{prop}}
\def\bo{\begin{observe}}
\def\eo{\end{observe}}
\def\bas{\begin{assumption}}
\def\eas{\end{assumption}}
\def\RR{\mathbb{R}}

\def\ZZ{\mathbb{Z}}
\def\NN{\mathbb{N}}

\def\beg{\begin{eg}}
\def\eeg{\end{eg}}

\def\p{\thinspace}

\makeindex

\numberwithin{equation}{section}
\numberwithin{table}{section}

\startlocaldefs
\endlocaldefs

\begin{document}

\begin{frontmatter}

\title{Limiting speed of a second class particle in ASEP}
\runtitle{Speed of a second class particle.}
\runauthor{Ghosal, Saenz \& Zell}

\begin{aug}
  


 \author{\fnms{Promit}  \snm{Ghosal}\textsuperscript{1}\corref{}\ead[label=e1]{pg2475@columbia.edu}\hspace{0.1in}}
  \author{\fnms{Axel} \snm{Saenz}\textsuperscript{2}\ead[label=e2]{ais6a@virginia.edu}\hspace{0.1in}}
  \author{\fnms{Ethan C.} \snm{Zell}\textsuperscript{3}\ead[label=e3]{ethanzell@virginia.edu}\vspace{0.1in}}

\address{\textsuperscript{1}Columbia University, Department of Statistics, 1255 Amsterdam, New York, NY 10027, USA\\ \printead{e1}}
\address{\textsuperscript{2}University of Virginia, Department of Mathematics, 141 Cabell Dr., Charlottesville, VA 22903, USA\\ \printead{e2}}
\address{\textsuperscript{3}University of Virginia, Department of Mathematics, 141 Cabell Dr., Charlottesville, VA 22903, USA\\ \printead{e3}}


\end{aug}

\begin{abstract}
     We study the asymptotic speed of a second class particle in the two-species asymmetric simple exclusion process (ASEP) on $\ZZ$ with each particle belonging either to the first class or the second class. For any fixed non-negative integer $L$, we consider the two-species ASEP started from the initial data with all the sites of $\ZZ_{<-L}$ occupied by first class particles, all the sites of $\ZZ_{[-L,0]}$ occupied by second class particles, and the rest of the sites of $\ZZ$ unoccupied. With these initial conditions, we show that the speed of the leftmost second class particle converges weakly to a distribution supported on a symmetric compact interval $\Gamma\subset \RR$. Furthermore, the limiting distribution is shown to have the same law as the minimum of $L+1$ independent random samples drawn uniformly from the interval $\Gamma$.    
\end{abstract}



\begin{keyword}
\kwd{Simple exclusion process}
\kwd{multispecies ASEP}
\kwd{block probability}
\end{keyword}

\end{frontmatter}









\tableofcontents

\section{Introduction}\label{introduction}
The goal of this paper is to study the asymptotic speed of a second class particle in the two-species asymmetric simple exclusion process (ASEP). The ASEP is an interacting particle system on $\ZZ$ with each of the vertices of $\ZZ$ occupied by at most one particle. Every particle carries an exponential clock of rate $1$, and all clocks are independent. The evolution of the particles is Markovian and may be described as follows: when the clock rings, the particle  decides to move to the right (or left) by one with probability $p\in (0,1)$ (or $1-p$). The particle moves to the corresponding target site if it was unoccupied, and nothing happens if the target site is occupied (i.e.~the jump is suppressed). In a two-species ASEP, every particle is labeled as either first class or second class, depending on its jump hierarchy. For instance, a first class particle can interchange its position with a second class particle which was sitting on a neighboring target site; on the other hand, a second class particle is not allowed to interchange its position with any other first or second class particle.

Fix $L\in \ZZ_{\geq 0}$. Consider the two-species ASEP with the right  jump probability $p\in (\frac{1}{2},1]$ and the left jump probability $q:=1-p$. Additionally, take infinitely many first class particles and $L+1$ second class particles in the system. Let $x_{n}(t)$ (resp.~$x^{*}_{n}(t)$) be the position of the first (resp.~second) class particle of index $n\in \ZZ_{\geq 0}$ at time $t\in \RR_{\geq 0}$. The initial positions of the particles are given by 
\begin{align}\label{eq:InitData}
x_{n} (0) &=
-n-L, \qquad  n \in \ZZ_{\geq 1}, \\
x^{*}_{n}(0) & = -n, \qquad\qquad  n\in \ZZ_{[0,L]}.    
\end{align}
 In the following, we state the main result of our paper.
 
\bt\label{thm:Main} 
 Given the asymmetry parameter $\gamma := p-q$ for $p \in (\frac{1}{2}, 1]$ and $q = 1-p$, we have  
\begin{align}\label{eq:ConvD}
\frac{x^{*}_{L}(t)}{t} \stackrel{d}{\to} U_{L}, \quad \text{as } t\to \infty 
\end{align}
with $U_L$, a random variable supported on $[-\gamma, \gamma]$, and 
\begin{align}\label{eq:Ul}
\mathbb{P}\big(U_L\geq s\big) = \Big(\frac{1-\gamma^{-1}s}{2}\Big)^{L+1}, \qquad \forall s\in [-\gamma, \gamma].  
\end{align} 
\et

Since its introduction (see \cite{MGP68}, \cite{Spitzer70}), the ASEP remains one of the most important interacting particle systems and serves as a testing ground for several conjectures about the KPZ universality class. We refer to \cite{Liggett85, Liggett99} for detailed historical accounts and further references. 

 Moreover, it is well known that that the hydrodynamic limit of the ASEP is governed by the \emph{invicid Burgers equation} which is written as 
 \begin{align}\label{eq:IBE}
 \frac{\partial u_t(x)}{\partial t} + \frac{\partial(u_t(x)(1-u_t(x))) }{\partial x} =0.
 \end{align}
  The characteristic line of this nonlinear PDE is intimately related with the trajectory of a second class particle. It was shown by many authors (for instance, \cite{Ferrari92,Reza91})       
that the second class particle started from a position $x$ follows the the characteristic line of \eqref{eq:IBE} started from $x$. If $x$ does not belong to the rarefaction fan, then there is only one characteristic line emanating from $x$. On the contrary, if the initial position is inside the rarefaction fan, then there are infinitely many characteristic lines coming out from $x$.  Theorem~\ref{thm:Main} essentially shows that the leftmost second class particle in the two-species ASEP started from \eqref{eq:InitData} will choose one of the characteristic lines\footnote{Each point of the interval $[-(p-q), (p-q)]$ denotes the slope of one of the characteristic lines.} from the interval $\Gamma:=[-(p-q),(p-q)]$ according to the law of the minimum of $L+1$ many random uniform samples from the interval $\Gamma$.

Study of the second class particle in the exclusion processes plays a unique role in its rich history. \cite{Reza95} obtained a law of large numbers for the position of the second class particle in the ASEP under \emph{shock initial data}. It is worth mentioning that the initial data of \eqref{eq:InitData} does not fall under the  category of the shock initial data. In fact, as time increases, we observe the formation of a rarefaction region where the density of the (first class) particles decreases linearly from $1$ to $0$ the initial data of \eqref{eq:InitData}. In the case of rarefaction fan, \cite{FK95} proved the weak convergence of the asymptotic speed of the second class particle (in the totally asymmetric simple exclusion process (TASEP)) to the uniform law on $[-1,1]$. The initial data of \cite{FK95} has a second class particle at $0$ and all the first class particles fully packed on one side of the origin. Setting $L=0$ and $p=1$ in Theorem~\ref{thm:Main} retrieves the main result of \cite{FK95}. Later, the result of \cite{FK95} is improved in \cite{MG05}; the authors showed the almost sure convergence of the asymptotic speed. See \cite{G14} for a related result in the totally asymmetric constant rate zero-range process. \cite{CD06} studied the Hammersely procss and obtained an exciting limit theorem for the second class particle at the rarefaction fan. \cite{FGM09} proved the analogue of the main result of \cite{FK95} in ASEP. They also studied the collision probability of two second class particles in the ASEP. \cite{TW09S} found the exact distribution of the second class particle in the ASEP. \cite{BN17} generalized the result of \cite{FK95} for a wide range of partially and totally asymmetric interacting particles systems.

Exploring the connection between the second class particle in the TASEP and the competition interface of the last passage percolation, many new results were found for the second class particle (see \cite{FP05, FMP09}). Building on this connection, \cite{CP13} derived a general formula for the limiting law of the asymptotic speed of the second class particle in the TASEP under arbitrary initial data. In particular, when the initial data of the TASEP has one second class particle at site $-L$ and all the other sites of $\ZZ_{\leq 0}$ are occupied by first class particles, it is a straightforward exercise to verify (using the formula of \cite{CP13}) that the limiting speed of the second class particle will be the same as \eqref{eq:Ul} (with $\gamma=1$). This brings us to the following conjecture which is supported by an extensive simulation shown in Appendix~. 
\begin{conj}
Fix $p\in (\frac{1}{2},1]$. Consider the ASEP on $\ZZ$ with the right jump probability $p$ started from the initial data with all the sites of $\ZZ_{\leq 0}\backslash \{0\}$ occupied by first class particles, except for a single second class particle at site $-L$, and the remaining sites unoccupied. Then, for some parameter $\alpha := \alpha(p,L)$ depending on $p$ and $L$ with $\alpha(1,L) =1$ and $\alpha(p, 0) = \gamma$, we have
\begin{align}
\frac{x^{*}_{L}(t)}{t} \stackrel{d}{\to} \tilde{U}_{L}, \quad \text{as } t\to \infty 
\end{align}
with $\tilde{U}_L$, a random variable supported on $[-\alpha, \alpha]$, and 
\begin{align}
\mathbb{P}\big(\tilde{U}_L\geq s\big) = \Big(\frac{1-\alpha^{-1}s}{2}\Big)^{L+1}, \qquad \forall s\in [-\alpha, \alpha].  
\end{align} 
\label{conjecture1}
\end{conj}
For further applications of the connection between the second class particle and the competition interface, we refer to \cite{RS,FNP17}. The TASEP with particles of multiple species (e.g., first class, second class, third class etc.) was studied in \cite{AAV11} where the authors found the joint limiting law of the speeds of the particles of different species. They also conjectured similar results for the multi-species ASEP. It is worth noting that Theorem~\ref{thm:Main} can be derived in an alternative way by assuming that Conjecture~1.9 and~1.10 of \cite{AAV11} hold.

There are three main tools used in our proof of Theorem~\ref{thm:Main}. Namely, $(a)$ a coupling between the two-species ASEP and the multi-species (or, colored) ASEP, $(b)$ a relation between the transition probabilities of the multi-species ASEP and the single-species (or, uncolored) ASEP\footnote{single-species ASEP contains only one class of particles.}, and $(c)$ an estimate of the $(L+1)$-block probability\footnote{Probability of having a block of $L+1$ consecutive particles.} of the single-species ASEP. In Lemma~\ref{lem:Coupling} of Section~\ref{sec:Multi}, we construct a coupling between the two-species ASEP started from  \eqref{eq:InitData} and the multi-species ASEP started from the step initial data (see Definition~\ref{bd:mASEP}) such that the position of the leftmost second class particle of the two-species ASEP agrees with the position of the leftmost particles among the first $L+1$ species of the multispecies ASEP for all $t\geq 0$. This brings us to the use of the transition probabilities of the multispecies ASEP. Recently, \cite{BorWhe} studied the spectral theory of the colored stochastic vertex models. They showed a correspondence between the  probabilities of the colored and uncolored stochastic six vertex models which under appropriate limit gives a similar description between the multispecies and the single-species ASEP. Briefly, \cite[Theorem~12.3.5]{BorWhe} connects the tail probabilities of the minimum postion of a finite set of particles in the multispecies ASEP with the \emph{block probabilities} of the single-species ASEP. We refer to Proposition~\ref{ppn:BorWhe} of Section~\ref{sec:Multi} for more details. Thanks to the coupling lemma and \cite[Theorem~12.3.5]{BorWhe}, our proof of Theorem~\ref{thm:Main} boils down to estimating the block probabilities of the single-species ASEP (which is provided by a recent work of \cite{TW18}). We refer to Proposition~\ref{ppn:TracyWidom} for more details on this estimate. Combining all these tools, we complete the proof of Theorem~\ref{thm:Main} in Section~\ref{sec:ThmMain}.


\section{Multispecies (or, colored) ASEP}\label{sec:Multi}
In this section, we review the the multispecies (or, colored) ASEP and construct a coupling between the latter and the two-species ASEP started from  the initial data \eqref{eq:InitData}. Towards the end of this section, we will recall the correspondence between the multi-species ASEP and the single-species ASEP.     

\bd\label{bd:mASEP} 
The \emph{multispecies ASEP} is a system of colored interacting particles on the infinite one-dimensional integral lattice, with sites labeled by $\ZZ$. Every particle in the system has a different color prescribed to it, which is a positive integer. Each site may be occupied by at most one particle. At a time $t$, a
configuration of the multi-species ASEP is given by 
\begin{align}
\eta^{\mathrm{mASEP}}(t) = \{\ldots , \eta_{-1}(t), \eta_{0}(t), \eta_{1}(t), \ldots \}, \quad \eta_i(t) \in \NN 
\end{align}
with $\eta_i(t) =0$ indicating that site $i$ is unoccupied at time $t$ and $\eta_i(t) =j>0$ indicating that site $i$ is occupied by a particle of color $j$ at time $t$. We denote the position of the particle with color $n$ by $x^{\mathrm{mASEP}}_n(t)$. As in ASEP, each particle carries an exponential clock, and all the clocks in the system are independent. When a clock rings, the corresponding particle jumps to the right with probability $p$ and to the left with probability $q:=1-p$. The particle interchanges positions with another particle if and only if the occupation variable $\eta_i(t)$ of the corresponding target site is lower. Otherwise, nothing happens.

\ed

\bd[Step initial data] The \emph{step initial data} is the initial configuration $\eta^{\mathrm{mASEP}} (0)$ with 
\begin{align}\label{eq:StepInit}
\eta^{\mathrm{mASEP}}_n(0) = 
\begin{cases}
-n & n\leq  0, \\
0 & n>0. 
\end{cases}
\end{align}
\ed

\bl\label{lem:Coupling}
 Fix $p\in (\frac{1}{2}, 1]$. Consider the two-species ASEP\footnote{We often refer this as \emph{two-color ASEP} in the proof.} on $\ZZ$ with the right jump probability $p$ started with the initial data \eqref{eq:InitData}. Additionally, consider  the multi-species (or, multi-color) ASEP started with the step initial data (see \eqref{eq:StepInit}). Then , there exists a coupling between these two systems such that 
\begin{align}\label{eq:Coupling}
x^{*}_{L}(t) = \min\{x^{\mathrm{mASEP}}_{n}(t): 0\leq n\leq L\}
\end{align}
for all $t\geq 0$.

\el
\smallskip 

\begin{proof}
 We start with an informal description of the coupling. Initially, for $n\in \ZZ_{>L}$, we couple the clock of the particle at $x_{n-L}(0)$ from the two-color ASEP with that of the particle at $x^{\mathrm{mASEP}}_{n}(0)$ from the multi-color ASEP. Similarly, we couple the clock of the second class particle at $x^{*}_{n}(0)$ with the particle at $x^{\mathrm{mASEP}}_{n}(0)$ for all $n \in \ZZ_{[0,L]}$. If two particles of the multi-color ASEP interchange their positions at time $t$ and both of them were coupled to the particles of same class of the two-color ASEP, we interchange the \emph{coupling status} of those two particles. That is, if $x^{\mathrm{mASEP}}_{r}(t-)$ is interchanged with $x^{\mathrm{mASEP}}_{s}(t-)$ at time $t$ and the clocks of the particles at $x^{\mathrm{mASEP}}_{r}(t-)$ and $x^{\mathrm{mASEP}}_{s}(t-)$ were coupled with those at $x_{r^{\prime}}(t-)$ (or, $x^{*}_{r^{\prime}}(t-)$) and $x_{s^{\prime}}(t-)$ (or, $x^{*}_{s^{\prime}}(t-)$) respectively, then, at time $t$, we couple the particles at $x^{\mathrm{mASEP}}_{r}(t)$ (resp.$x^{\mathrm{mASEP}}_{s}(t)$) with $x_{s^{\prime}}(t)$ (resp. $x_{r^{\prime}}(t)$). If the particles of the multi-species ASEP participating in an interchange are coupled to particles from separate classes in the two-colored ASEP, then we interchange the corresponding particles of the two-color ASEP. As a consequence, the \emph{coupling status} of the particles does not change, i.e., the set of the pairs of particles which were coupled to each other before remains same. For further clarity, we illustrate this coupling with the following figure.
 
\begin{figure}[h]
\centering
\begin{subfigure}{.4\textwidth}
\centering
  \includegraphics[width=0.8\linewidth]{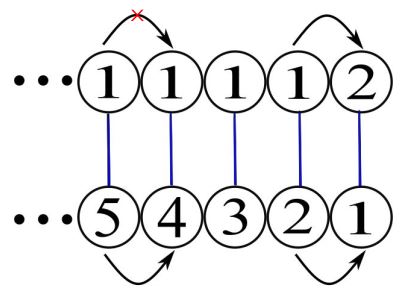}
  \label{fig:sub1}
\end{subfigure}
\begin{subfigure}{.4\textwidth}
  \centering
  \includegraphics[width=0.75\linewidth]{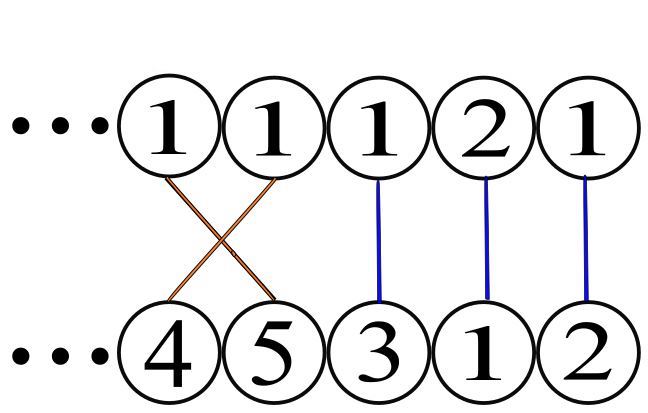}
  \label{fig:sub2}
\end{subfigure}
\caption{The left image represents the configuration at $t-$ and the right image at time $t$. In the right image, the first example of movement swaps coupling; in the second example, it does not.}
\label{fig:test}
\end{figure}
 

Now, we are ready to give an explicit description of the coupling.  
Set
\begin{align}
\tau_0 := \inf\{t\geq 0: x^{\mathrm{mASEP}}_{m}(t) = x^{\mathrm{mASEP}}_{m}(t-)+1, x^{\mathrm{mASEP}}_{m-1}(t)= x^{\mathrm{mASEP}}(t)-1 \text{ for some }m\in \ZZ_{\geq 1} \}. 
\end{align}
In other words, $\tau_0$ is the first time when two particles of the multi-species ASEP interchange their positions. Then, for all $t<\tau_0$, we couple the exponential clocks of the particles of the two-color ASEP and the multi-color ASEP in the following manner:   
 \begin{align}
 x_{n-L}(t)  =
 x_{n-L}(t-)\pm 1, &\quad  \text{ if } x^{\mathrm{mASEP}}_{n}(t) = x^{\mathrm{mASEP}}_{n}(t-) \pm 1 \text{ and } n\in \ZZ_{>L},\label{eq:InitCoup1} \\
 x^{*}_{n}(t)  =
 x^{*}_{n}(t-)\pm 1,
  &\quad \text{ if } x^{\mathrm{mASEP}}_{n}(t) = x^{\mathrm{mASEP}}_{n}(t-) \pm 1  \text{ and } n\in \ZZ_{[0,L]}.\label{eq:InitCoup2}
 \end{align}

 To encode the coupling between the particles of the two-color ASEP and the particles of multi-color ASEP, we introduce a random function $f:[0,\infty)\times \ZZ \to \NN\cup\{0^{*}, 1^{*}, \ldots , L^{*}\}$\footnote{Here, $n^{*}$ stands for the index of the second class particle of the site $x^{*}_n(0)$ at time $t=0$.} which we often refer to as the \emph{coupling status}. In other words, $(n,f_t(n))$ denotes the pair of indices of the particles from the two-color ASEP and the multi-color ASEP which are coupled with each other at time $t$. For all $t<\tau_0$, we define 
 \begin{align}
 f_t(n) = \begin{cases}
 n-L & n \in \ZZ_{>L},\\
 n^{*} & n \in \ZZ_{[0,L]}.
 \end{cases}
 \end{align}

  Now, we extend the coupling for all $t\geq \tau_0$. Consider a sequence of stopping times $\{\tau_1< \tau_2< \ldots \}$ given by
\begin{align}
\tau_{k} := \inf\big\{t>\tau_{k-1}: x^{\mathrm{mASEP}}_{r}(t) = x^{\mathrm{mASEP}}_{s}(t-), x^{\mathrm{mASEP}}_{s}(t) = x^{\mathrm{mASEP}}_{r}(t-) \text{ for some }r> s\geq 1\big\}.
\end{align}
Due to the independence between the exponential clocks associated to the particles, the probability of interchanging the positions of two pairs of particles at the same time is $0$. We update $f_t$ for all $t=\tau_0, \tau_1, \tau_2, \ldots $ and set it to be constant for all $t\in [\tau_k, \tau_{k+1})$ and for all $k\geq 0$. If the position of the particles of color $r$ and $s$ are interchanged for some $r>s$ at $t=\tau_k$, we define
\begin{align}
f_{\tau_k}(n) = \begin{cases}
f_{\tau_k-}(n) & n \neq r,s,\\
f_{\tau_k-}(s)\mathbbm{1}(\mathcal{A}_{r,s}) +f_{\tau_k-}(r)\mathbbm{1}(\mathcal{A}^c_{r,s}) & n=r, \\ 
f_{\tau_k-}(r)\mathbbm{1}(\mathcal{A}_{r,s}) +f_{\tau_k-}(s)\mathbbm{1}(\mathcal{A}^c_{r,s}) & n=s. 
\end{cases} \label{eq:InChange}
\end{align}
with
\begin{align}
\mathcal{A}_{r,s} &:= \Big\{(f_{\tau_k-}(r),f_{\tau_k-}(r))=(k,\ell), \text{ for some }k, \ell\in \ZZ_{\geq 1}\Big\}\\&\cup \Big\{(f_{\tau_k-}(r),f_{\tau_k-}(r))=(k^{*},\ell^{*}), \text{ for some }k, \ell\in \ZZ_{[0,L]}\Big\}.
\end{align}
For any $t\in (\tau_k, \tau_{k+1})$ and $k\geq 0$, we couple the two-colored ASEP and the multi-colored ASEP in the following way: 
\begin{align}\label{eq:Basic}
x_{f_{t}(n)}(t) = x_{f_{t}(n)}(t-)\pm 1, \quad\text{if}\quad   x^{\mathrm{mASEP}}_{n}(t) = x^{\mathrm{mASEP}}_{n}(t-)\pm 1. 
\end{align} 
This completes the description of the coupling.\par

 Now, we show that under this coupling, \eqref{eq:Coupling} holds for all $t\geq 0$. Note that, by \eqref{eq:InitCoup1}-\eqref{eq:InitCoup2}, we immediately have
$$  \min\{x^{\mathrm{mASEP}}_{n}(t):1\leq n\leq \ell\} = x^{\mathrm{mASEP}}_{\ell}(t)=x^{*}_{\ell}(t), \quad \forall \ell\in \ZZ_{[0,L]}$$  
for all $t< \tau_0$. This shows that \eqref{eq:Coupling} holds true if $t$ is less than $\tau_0$. Since the set of coupled pairs $\{(n, f_t(n))\}_{n\geq 1}$ remains unchanged for all $t\in (\tau_k, \tau_{k+1})$, it suffices to show that \eqref{eq:Coupling} holds at $t=\tau_k$ for all $k\geq 0$. We show this by induction. 

First, we prove the claim for $t=\tau_0$. Suppose two particles of color $r$ and $s$ interchange their positions at $t=\tau_0$. If neither $f_{\tau_0-}(r)$ nor $f_{\tau_0-}(s)$ belongs to $ \{0^*, 1^{*}, \ldots , L^{*}\}$, then \eqref{eq:Coupling} holds at $t=\tau_{0}$ because the \emph{coupling status} of the second class particles does not change at $t=\tau_0$. If $f_{\tau_0-}(r)\in \{0^*, 1^{*}, \ldots , L^{*}\}$ but $f_{\tau_{0}-}(r)\notin \{0^*, 1^{*}, \ldots , L^{*}\}$ (or vice versa), then again \eqref{eq:Coupling} holds at $t=\tau_0$ for the same reason. Lastly, consider the case when both $f_{\tau_0-}(r)$ and  $f_{\tau_0-}(s)$ belong to $\{0^*, 1^{*}, \ldots , L^{*}\}$. Suppose $f_{\tau_0-}(r) = L^{*}$ and $f_{\tau_0-}(s)= \ell^{*}$ for some $\ell\in \ZZ_{[0, \ell)}$. Since \eqref{eq:Coupling} holds for $t<\tau_0$, we have $x^{\mathrm{mASEP}}_{r}(\tau_{0}-)<x^{\mathrm{mASEP}}_{s}(\tau_{0}-)$. Appealing to \eqref{eq:InChange}, we observe  
\begin{align}
x^{*}_{L}(\tau_0) = x^{\mathrm{mASEP}}_{s}(\tau_0)& =  x^{\mathrm{mASEP}}_{r}(\tau_0-)\\&=\min\{x^{\mathrm{mASEP}}_{n}(\tau_0-):0\leq n\leq L\} =\min\{x^{\mathrm{mASEP}}_{n}(\tau_0):0\leq n\leq L\} 
\end{align}
with the last equality following from the fact that the minimum of a set of indexed integers is invariant under the mutual interchange of the indices. This verifies that \eqref{eq:Coupling} holds at $t=\tau_0$.

  Let us now assume \eqref{eq:Coupling} holds for all $t=\tau_0, \ldots , \tau_{k_0}$. We need to show that it will also hold for $t=\tau_{k_0+1}$. However, the proof of this claim is almost verbatim to the proof of the same claim for $t=\tau_0$, which we do not repeat. This completes the proof.  

\end{proof}

\begin{rem}
It is worth mentioning that one can construct a coupling between the two-species TASEP and the multi-species TASEP in a similar way as in Lemma~\ref{lem:Coupling} such that \eqref{eq:Coupling} holds even when the initial data for the two-species TASEP has only one second class particle at $-L$ and the rest of the sites of $\ZZ_{\leq 0}$ are occupied by the first class particles. Indeed, using this coupling and similar arguments as in Section~\ref{sec:ThmMain}, one can prove \eqref{eq:ConvD} for two-species TASEP started from this special initial data.      
\end{rem}

We end this section with the correspondence between the joint probability of the multi-species ASEP and the block probability of the single-species ASEP. 

Consider two sets of integers $\mathcal{I}:=\big\{I_1<\ldots < I_{k}\leq P\big\}$ and  $\mathcal{J}:=\big\{1\leq J_1< \ldots < J_{\ell}\big\}$ for some $\ell\in \NN$ and $P_0\in \ZZ$. Moreover, consider 
\begin{align}
\mathbb{P}^{\mathrm{asep}} (\mathcal{I}, \mathcal{J}; P,t) := \mathrm{Prob}\Big(\eta^{\mathrm{asep}}_{I_i}(t) =1, \forall 1\leq i\leq k, \text{ and }\eta^{\mathrm{asep}}_{P+J_{j}}(t)=1, \forall 1\leq j\leq \ell\Big),
\end{align}
which is the probability that, in the single-species ASEP (started from the \emph{step initial data}\footnote{The step initial data $\eta^{\mathrm{asep}}$ for the single-species ASEP is given by $\eta^{\mathrm{asep}}_{n}(0)=1$ if $n\in \ZZ_{\leq 0}$ and $\eta^{\mathrm{asep}}_{n}(0)=0$ if $n\in \ZZ_{>0}$.}) at time $t$, each of the positions $\big\{I_1, \ldots , I_k\big\}$ and $\big\{P+J_1, \ldots , P+J_{\ell}\big\}$ is occupied by a particle; and 
\begin{align}
\mathbb{P}^{\mathrm{mASEP}}(\mathcal{I}, \mathcal{J}; P,t) := \mathrm{Prob}\Big(\eta^{\mathrm{mASEP}}_{I_i}(t)\geq 1, \forall 1\leq i\leq k,\p\&\p\exists p_j>P\text{ s.t. }\eta^{\mathrm{mASEP}}_{p_j}(t)= J_j, \forall 1\leq j\leq \ell\Big)
\end{align} 
describing the probability that, in the multi-species ASEP (started from the step initial data) at time $t$, each of the positions $\{I_1, \ldots, I_k\}$
is occupied by a particle and the particles of colors $\{J_1,\ldots , J_{\ell}\}$ are all situated strictly to the right of position $P$.

\bp[Theorem~12.3.5 of \cite{BorWhe}]\label{ppn:BorWhe}
For all $t\geq 0$, any $P_0\in \ZZ$ and $\mathcal{I}=\big\{I_1<\ldots < I_{k}\leq P_0\big\}$, $\mathcal{J}:=\big\{1\leq J_1< \ldots < J_{\ell}\big\}$, 
\begin{align}\label{eq:coupling}
\mathbb{P}^{\mathrm{asep}} (\mathcal{I}, \mathcal{J}; P,t) = \mathbb{P}^{\mathrm{mASEP}}(\mathcal{I}, \mathcal{J}; P,t).
\end{align} 
\ep

\section{Proof of Theorem~\ref{thm:Main}}\label{sec:ThmMain}

%

We combine \eqref{eq:coupling} with the coupling of Lemma~\ref{lem:Coupling} to prove our main result. Fix  $\mathcal{J}:= \{1, 2,3, \ldots , L+1\}$ and $\mathcal{I} = \{-k\}$ for some $k>L$. From \eqref{eq:Coupling} in Lemma~\ref{lem:Coupling}, we have
\begin{align}
\Big\{ &x^{*}_{L}(t) >  s t, \p\p \eta_{-k}(t) =1\Big\}\\ &= \Big\{ x^{\mathrm{mASEP}}_{j}(t)> st, \p\p \forall 0\leq j\leq L, \quad \eta^{{\mathrm{mASEP}}}_{-k}(t) \geq 1\Big\} \\
& = \Big\{\eta^{{\mathrm{mASEP}}}_{k}(t)\geq 1, \exists p_j > s t \p\p\text{ s.t } \eta^{{\mathrm{mASEP}}}_{p_j}(t) = j, \forall 0\leq j\leq L \Big\}\label{eq:Equiv}
\end{align}
Also, due to \eqref{eq:Equiv} and \eqref{eq:coupling}, we have
\begin{align}\label{eq:Duality}
\mathbb{P}^{\mathrm{ASEP}}&\Big(\Big\{x^{*}_{L}(t) \geq  s t,\p\p \eta_{-k}(t) =1\Big\}\Big)\\& =\mathbb{P}^{\mathrm{mASEP}}\Big(\Big\{\eta^{{\mathrm{mASEP}}}_{k}(t)\geq 1, \exists p_j >s t \text{ s.t } \eta^{{\mathrm{mASEP}}}_{p_j}(t) = j, \forall 0\leq j\leq L \Big\}\Big)\\& = \mathbb{P}^{\mathrm{asep}}\Big(\eta^{\mathrm{asep}}_{k}(t) =1 , \text{ and }\eta^{\mathrm{asep}}_{\sigma t+j}=1, \forall 0\leq j\leq L\Big).  
\end{align}
Here, $\mathbb{P}^{\mathrm{ASEP}}$ denotes the probability with respect the two-species ASEP started from \eqref{eq:InitData} whereas $\mathbb{P}^{\mathrm{asep}}$ stands for the probability with respect to the uncolored ASEP started from the step initial condition. Letting $k$ go to $\infty$ on both sides of \eqref{eq:Duality}, we obtain 
\begin{align}\label{eq:IntStep}
\mathbb{P}^{\mathrm{ASEP}}\Big(x^{*}_{L}(t) \geq  s t\Big) = \mathbb{P}^{\mathrm{asep}}\Big(\eta_{s t+j}=1, \forall 0\leq j\leq L\Big). 
\end{align}
To complete the proof of the main result, we introduce the following lemma. Note that combining \eqref{eq:RLim} with \eqref{eq:IntStep} completes the proof of Theorem~\ref{thm:Main}.
\smallskip 

\bl  \label{lm:limit_block_prob}
Given the asymmetry parameter $\gamma : = p-q$ with $p \in (\frac{1}{2}, 1]$ and $q =1-p$, we have
\begin{align}\label{eq:RLim}
\lim_{t\to \infty}\mathbb{P}^{\mathrm{asep}}\Big(\eta_{s t+j}=1, \forall 0\leq j\leq L\Big) = \Big(\frac{1-\gamma^{-1}s}{2}\Big)^{L+1}
\end{align}
for any $L\geq 0$.
\el
\smallskip 

For proving Lemma~\ref{lm:limit_block_prob}, we need the following result from \cite{TW18}.

\bp[Theorem~1 of \cite{TW18}]\label{ppn:TracyWidom} Consider single-species ASEP on $\ZZ$ (with the right jump probability $p\in (\frac{1}{2},1)$) started from the step initial data. Denote the position of the particle of index $n\in \ZZ_{\geq 0}$ at time $t\in \RR_{\geq 0}$ by $x^{\mathrm{asep}}_{n}(t)$. Let $\mathcal{P}_{L, \mathrm{step}}(x,m,t)$ be the probability of the event 
\begin{align}
x^{\mathrm{asep}}_{m}(t) = x,\quad  x^{\mathrm{asep}}_{m+2}(t)=x+1, \quad x^{\mathrm{asep}}_{m+2}(t) = x+2, \quad \dots \quad  ,\quad x^{\mathrm{asep}}_{m+L}(t)= x+L. 
\end{align} 
Additionally, we set $m = \sigma t$ for some $\sigma \in (0,1)$ and introduce the parameters
	\begin{equation}\label{eq:parameters}
	c_1 = 1 - 2 \sqrt{\sigma} \quad \text{and} \quad c_2 = \sigma^{-1/6}(1 - \sqrt{\sigma})^{2/3}.
	\end{equation}
Then, for $x = c_1 t - c_2 \zeta t^{1/3}$, one has 
\begin{align}\label{eq:BlockProbability}
\mathcal{P}_{L, \mathrm{step}}\Big(x, m, t/ \gamma \Big) = c^{-1}_2 \sigma ^{(L-1)/2} F^{\prime}_{\mathrm{GUE}}(\zeta) t^{-\frac{1}{3}} + o(t^{-\frac{1}{3}})
\end{align}
with $F^{\prime}_{\mathrm{GUE}}$, the derivative of the Tracy-Widom GUE distribution\footnote{It is the limiting distribution of the largest eigenvalue of the Gaussian unitary ensemble.}.
\ep

\begin{proof} [Proof of Lemma \ref{lm:limit_block_prob}]
Let us define $\theta$ and $\zeta$ implicitly by the following equalities:
 \begin{align}\label{eq:Notation}
\frac{s}{\gamma}=: c_1 + \theta c_2 (\gamma t)^{-\frac{2}{3}}, \p\p c_1:= 1-2\sqrt{\sigma+ \zeta (\gamma t)^{-\frac{2}{3}}}, \p\p \sigma :=\big(\frac{1-s/\gamma}{2}\big)^2, \p\p  c_2:= \sigma^{-\frac{1}{6}} (1-\sqrt{\sigma})^{\frac{2}{3}},
 \end{align}
with parameters $c_1$, $c_2$, and $\sigma$ as in \eqref{eq:parameters} in Proposition \ref{ppn:TracyWidom}. Expanding $c_1$, we have
\begin{align}
\frac{s}{\gamma}=c_1+\theta c_2 (\gamma t)^{-\frac{2}{3}} = 1-2\sqrt{\sigma} - (\sigma)^{-\frac{1}{2}}\zeta (\gamma t)^{-\frac{2}{3}} + \theta c_2 (\gamma t)^{-\frac{2}{3}} +o(t^{-\frac{2}{3}}). 
\end{align} 
Then, we deduce $\theta =\theta(\zeta)= (\sqrt{\sigma}c_2)^{-1}\zeta + o(1)$. Now, we fix $\zeta>0$ large. Also, we introduce some short-hand notations:
\begin{align}
\mathbb{P}_1 &: =\mathbb{P}^{\mathrm{asep}} \Big(\eta^{\mathrm{asep}}_{st+j}(t)=1, \forall 1\leq j\leq L\Big)\\
\mathbb{P}_2 & : = \mathbb{P}^{\mathrm{asep}}\Big(\Big\{\eta^{\mathrm{asep}}_{st+j}(t)=1, \forall 1\leq j\leq L \Big\}\cap \Big\{x^{\mathrm{asep}}_{\sigma \gamma t + \zeta (\gamma t)^{\frac{1}{3}}}(t)<st\Big\}\cap \Big\{x^{\mathrm{asep}}_{\sigma\gamma t - \zeta (\gamma t)^{\frac{1}{3}}}(t)>st+L\Big\}\Big)\\
\mathbb{P}_3 & := \mathbb{P}^{\mathrm{asep}}\Big(x^{\mathrm{asep}}_{\sigma \gamma t + \zeta (\gamma t)^{\frac{1}{3}}}(t)\geq st\Big), \qquad \mathbb{P}_4  := \mathbb{P}^{\mathrm{asep}}\Big(x^{\mathrm{asep}}_{\sigma\gamma t - \zeta (\gamma t)^{\frac{1}{3}}}(t)\leq st+L\Big)
\end{align}
It straightforward to see that 
\begin{align}
|\mathbb{P}_1- \mathbb{P}_2|\leq \mathbb{P}_3 + \mathbb{P}_4. 
\end{align}  
Note that, by setting $s = \gamma(c_1 + \theta c_2 t^{-2/3})$, we have
\begin{align}
\mathbb{P}_3 &\leq \mathbb{P}^{\mathrm{asep}}\Big(x^{\mathrm{asep}}_{\sigma \gamma t + \zeta(\gamma t)^{\frac{1}{3}}}(t)\geq (1-2\sqrt{\sigma + \zeta (\gamma t)^{-\frac{2}{3}}})\gamma t +\theta c_2(\gamma t)^{\frac{1}{3}} \Big)\label{eq:P3Bd}\\
\mathbb{P}_4 &\leq \mathbb{P}^{\mathrm{asep}}\Big(x^{\mathrm{asep}}_{\sigma \gamma t -\zeta(\gamma t)^{\frac{1}{3}}}(t)\leq \big(1-2\sqrt{\sigma - \zeta(\gamma t)^{-\frac{2}{3}}})\gamma t -\theta c_2 (\gamma t)^{\frac{1}{3}} \Big)\label{eq:P4Bd}
\end{align}
Moreover, by the result in \cite[Theorem~3]{Tracy2009b}, we know 
\begin{align}
\lim_{t\to \infty} \text{r.h.s of \eqref{eq:P3Bd}} = 1-F_{\mathrm{GUE}}(\theta), \qquad \lim_{t\to \infty} \text{r.h.s of \eqref{eq:P3Bd}} = F_{\mathrm{GUE}}(-\theta)
\end{align} 
with $F_{\mathrm{GUE}}$ denoting the Tracy-Widom GUE distribution. Therefore, 
\begin{align}
\lim_{\zeta\to \infty} \lim_{t\to \infty}\max\{\mathbb{P}_3, \mathbb{P}_4\} = 0
\end{align} 
since $\theta$ increases to $\infty$ as $\zeta\to \infty$. This implies 
\begin{align}
\lim_{t\to \infty} \mathbb{P}_1 = \lim_{\zeta\to \infty}\lim_{t\to \infty} \mathbb{P}_2. 
\end{align}
Now, take $\mathcal{P}_{L, \mathrm{step}}(x,m,t)$ (from Proposition \ref{ppn:TracyWidom}), the probability of the event 
\begin{align}
x^{\mathrm{asep}}_{m}(t) = x,\quad  x^{\mathrm{asep}}_{m+2}(t)=x+1, \quad x^{\mathrm{asep}}_{m+2}(t) = x+2, \quad \dots \quad, \quad x^{\mathrm{asep}}_{m+L}(t) = x+L. 
\end{align}
Then, we write 
 \begin{align}
 \mathbb{P}_2 = \sum_{j=\lfloor \sigma (\gamma t) - \zeta(\gamma t)^{\frac{1}{3}}\rfloor}^{\lceil \sigma (\gamma t) + \zeta(\gamma t)^{\frac{1}{3}}\rceil}\mathcal{P}_{L, \mathrm{step}}(st,j,t). 
 \end{align}
We approximate the sum above by the integral
\begin{align}
\mathbb{P}_2 = \sigma^{\frac{L+1}{2}}\int_{-c^{-1}_2\sigma^{-1/2}\zeta}^{c^{-1}_2\sigma^{-1/2}\zeta} F^{\prime}_{\mathrm{GUE}}(\xi)d\xi + o(1)
\end{align}
with $j = \lfloor \sigma\gamma t + c_2\sigma^{\frac{1}{2}}\xi (\gamma t)^{\frac{1}{3}}\rfloor$ for some $\xi \in [-\zeta, \zeta]$ and using \eqref{eq:BlockProbability} from Proposition \ref{ppn:TracyWidom}. Letting $t\to \infty$ and $\zeta\to \infty$ on both sides, we obtain
\begin{align}
\lim_{\zeta\to \infty}\lim_{t\to \infty} \mathbb{P}_2 = \sigma^{(L+1)/2}\int^{\infty}_{-\infty} F^{\prime}_{\mathrm{GUE}}(s) ds = \sigma^{(L+1)/2} 
\end{align}
with the last equality following by noting that $\int^{\infty}_{-\infty} F^{\prime}_{\mathrm{GUE}}(s) ds = F_{\mathrm{GUE}}(\infty) - F_{\mathrm{GUE}}(-\infty)= 1$. This shows \eqref{eq:RLim} and thus, completes the proof of the lemma. 
\smallskip

\end{proof}

\section{Appendix}\label{sec:Appendix}

In this section, we outline the reasoning responsible for the formulation of Conjecture \ref{conjecture1}. Namely, we utilized computer programmed simulations to model ASEP with a single second-class particle and used the simulations to estimate the cumulative distribution function (CDF) of the location of the second-class particle. We utilized \emph{Sage}, a well-known programming language akin to \emph{Python} but with many more mathematical libraries, to code the simulations; \emph{Rivanna}, the University of Virginia's High-Performance Computing (HPC) system, to run the simulations; and \emph{Mathematica} to analyze the data from the simulations.

In the case $p=1$ $(q=0),$ the simulations correspond to the TASEP with a second-class particle---the Sage program produced a list of ordered pairs from $10,000$ trials in order to generate an empirical CDF for the location of the second class particle after a fixed time $t=500$. In Mathematica, we analyzed the CDFs by taking numerical derivatives of the plots, and this allowed us to determine the degree of a polynomial to fit to the experimental data.  Our experimental results matched \cite{CP13} for the first starting positions of the second class particle. Figure \ref{fig:sub1} is one example with $x^*_2(0)=-2$. The experimental data (blue) matches the theoretical result (green), and the best quadratic fit to the experimental data is given in (red).

\begin{figure}[h]
\centering
\begin{subfigure}{.4\textwidth}
  \includegraphics[width=1.2\linewidth]{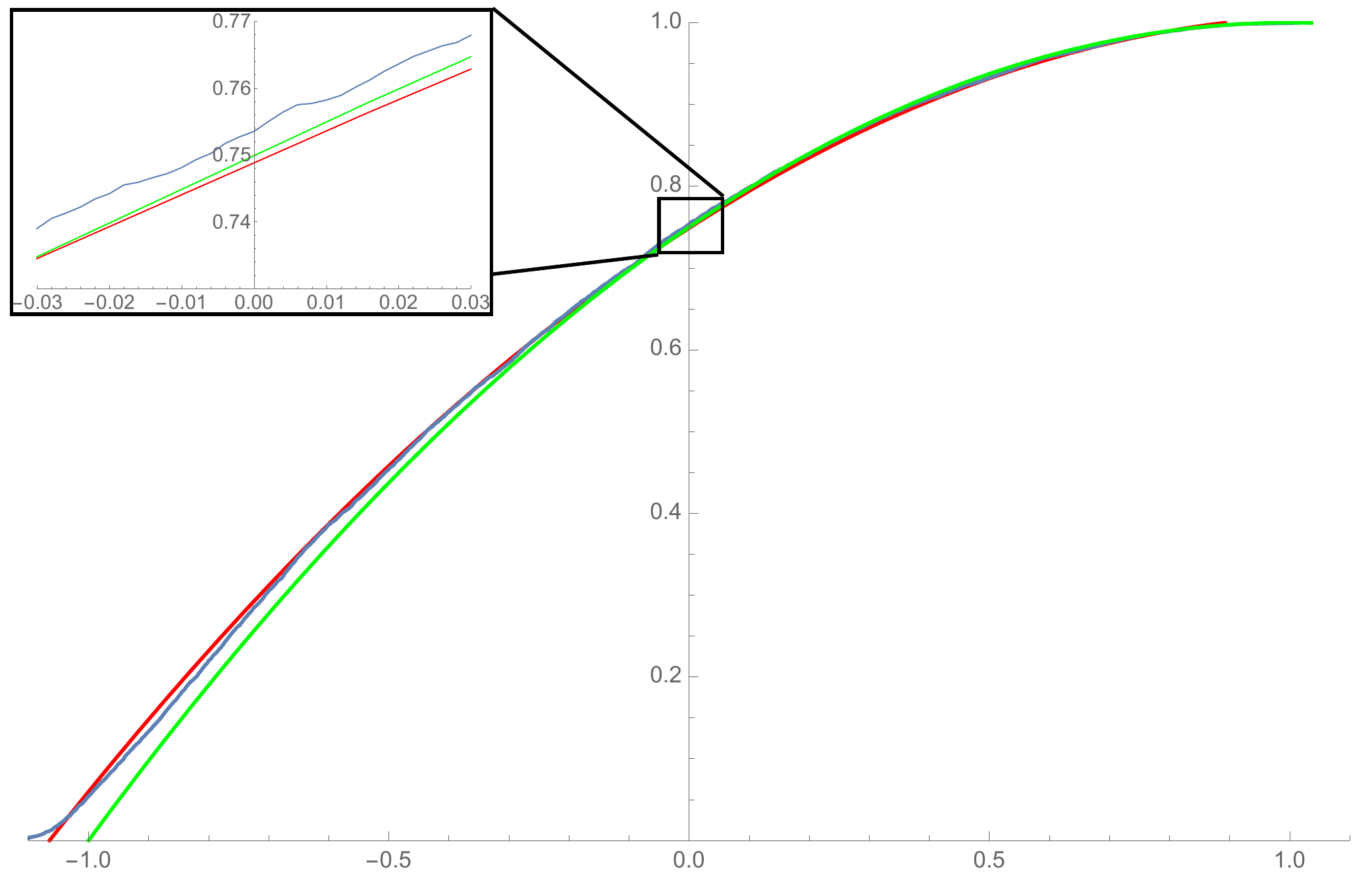}
  \caption{$p=1$}
  \label{fig:sub1}
\end{subfigure}
\hspace{5mm}
\begin{subfigure}{.4\textwidth}
  \centering
  \includegraphics[width=1.2\linewidth,trim=0 0 0 -1mm]{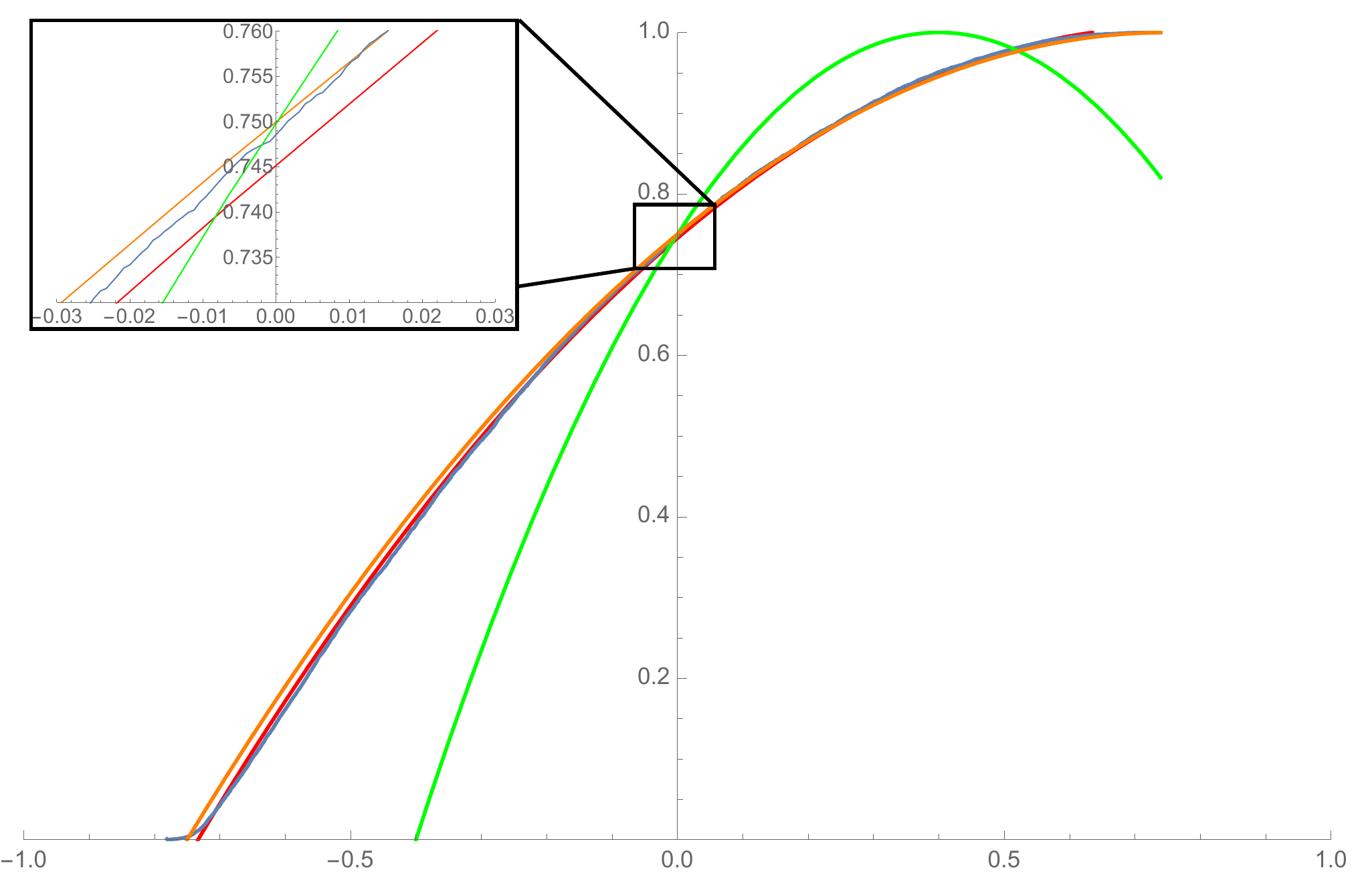}
  \caption{$p=0.7$}
  \label{fig:sub2}
\end{subfigure}
\caption{}
\label{fig:test}
\end{figure}

Hoping to detect a similar result for the ASEP with a single second class particle, we ran an identical analysis for various values of $p\in (\frac{1}{2},1).$ However, the experimental data for ASEP indicated that the TASEP result does not extend by simply introducing the asymmetry parameter $\gamma = p -q$ as in Theorem \ref{thm:Main}. We consider the following example. Let $p=0.7$ (and $q=0.3$), $\gamma = 0.4$, $x^*_2(0)=-2,$ and run at least $10,000$ trials. With the experimental data in blue, a quadratic best-fit in red, the polynomial $((1-\alpha^{-1}x)/2)^2$ in orange for $\alpha=\frac{7}{8},$ and the polynomial $((1-\gamma^{-1}x)/2)^2$ in green, we obtained \eqref{fig:sub2}.

Qualitatively, it is evident that extending the result for the TASEP by by simply introducing the asymmetry parameter $\gamma = p -q$, indicated in green, does not match the experimental data for ASEP. However, the orange line is promisingly close to the experimental data; namely, it is still possible that the CDF for the second-class particle in the ASEP case has the form $((1-\alpha^{-1}x)/2)^2,$ with $\alpha$ dependending on $p$ but given by a different formula other than $2p-1=\gamma.$ In the future, we hope to carry out simulations for many values of $p$ to elucidate a potential formula for $\alpha.$ 

\paragraph{Acknowledgments.}
We are grateful to Amol Aggarwal, Ivan Corwin, Patrick Ferrari, Leonid Petrov, and Craig Tracy for helpful discussions. This work was originally started as a research project for the \emph{Undergraduate Research Program} at the Mathematics Department in the University of Virginia. We are very grateful to the faculty organizers Julie Bergner and Thomas Mark. Also, we would like to thank the undergraduate participants Eric Keener, Cedric Harper, and Fernanda Yepez-Lopez for the tremendous help in learning the subject at the beginning stages of this project. Additionally, P.G. and A.S. would like to thank the organizer of the \emph{Integrable Probability Focused Research Group} (funded by NSF grants DMS-1664531, 1664617, 1664619, 1664650) for organizing stimulating events that propelled this project forward. A.S. was partially suppported by NSF grants DMS-1664617. E.C.Z was partially supported by NSF grant DMS-1659931.

 \bibliographystyle{plainnat}

\end{document}